\documentclass[reqno,12pt]{amsart}
\usepackage{amssymb,graphicx}
\def\R{{\mathbb R} }

\def\diam{{\rm diam\,}}

\newtheorem{Thm}{Theorem} 

\newtheorem{Prop}[Thm]{Proposition}

\newenvironment{Proof}[0]
{\medskip \noindent {\it Proof.} \ }{\ \hfill $\Box$\medskip}

\begin{document}
\title  {Local structure of self-affine sets}
\author{Christoph Bandt and Antti  K\"aenm\"aki}

\begin{abstract} The structure of a self-similar set with open set condition
does not change under magnification. For self-affine sets the situation is completely different.
We consider self-affine Cantor sets $E\subset\R^2$ of the type studied by Bedford, McMullen, Gatzouras and Lalley, for which the projection onto the horizontal axis is an interval.
We show that in small square $\varepsilon$-neighborhoods $N$ of almost each point $x$ in $E,$ with respect to many Bernoulli measures on address space, $E\cap N$ is well approximated by product sets $[0,1]\times C$ where $C$ is a Cantor set. 
Even though $E$ is totally disconnected,
all tangent sets have a product structure with interval fibres, reminiscent of the view of attractors of chaotic differentiable dynamical systems. We also prove that $E$ has uniformly scaling scenery in the sense of Furstenberg, Gavish and Hochman: the family of tangent sets is the same at almost all points $x.$
\end{abstract}
\maketitle
 \vspace{12mm}\\
Christoph Bandt\\
Institute for Mathematics and Informatics\\
Arndt University\\ 17487 Greifswald, Germany\\
e-mail: {bandt@uni-greifswald.de} \vspace{3mm}\\
Antti  K\"aenm\"aki\\
Department of Mathematics and Statistics\\
P.O. Box 35 (MaD)\\
40014 University of Jyv\"askyl\"a, Finland\\
e-mail: {antti.kaenmaki@jyu.fi}
\pagebreak

\section{Introduction}
The local fine structure of a fractal is intricate and regular at the same time.  In general, it is not hard to show that there are no tangent planes (cf. \cite{BK}, Theorem 1). Successive magnification of the set around a given point $x,$ the so-called {\em scenery flow} \cite{BFU,Ga,Hoc}, will never lead to a unique ``tangent set''. On the other hand, there is a well-defined {\em family of tangent sets,} termed  the ``tangential measure distribution'' \cite{Gr,MP} or ``gallery of micro-sets''  \cite{Fu} of the fractal. This family exists globally: it is obtained from the magnification flow at almost every point $x.$  This was proved for self-similar \cite{Gr,Ban} and self-conformal \cite{BFU} measures, including random constructions \cite{X,Fu}.  Recent papers by Furstenberg, Gavish, Hochman and Shmerkin \cite{Fu,Ga,Hoc} show that such a ``uniformly scaling scenery'' implies a number of nice geometric properties of the given fractal set or measure.

For self-similar sets with open set condition this is easy to understand since the structure does not change under magnification. 
Tangent sets are essentially the same as parts of the set. In Furstenberg's terminology, such fractals are homogeneous since all micro-sets are mini-sets \cite{Fu}. For self-conformal sets, the situation is only slightly different: tangent sets relate to parts of the fractal in a similar way as a tangent to a corresponding curve.

Here we study certain self-affine sets, and show in Theorem 2 that they  have a ``uniformly scaling scenery''. Our main result, however, is that they undergo a metamorphosis when they are magnified. Disconnected sets will turn into connected tangent sets. We show that for certain totally disconnected self-affine sets, the tangent sets have a product structure with connected fibres, like attractors of differentiable dynamical systems. The important assertion is not that such tangent sets exist -- this is fairly easy to see. The surprising fact is that essentially \emph{all tangent sets have this form.}
Actually, we conjecture that this is the typical local structure for large classes of self-affine sets,  with and without the open set condition.  \smallskip

\emph{Basic definitions.} We consider a family of plane self-affine Cantor sets in the unit square $Q=[0,1]^2$ which includes the types studied by Bedford \cite{Be}, McMullen \cite{MM}, Gatzouras and Lalley \cite{GL}. 
For $j\in J=\{ 1,...,m\},$ let $f_j:Q\to Q$ denote contractive map of the form
 \[  f_j(x_1,x_2)=(r_jx_1, s_jx_2)+(a_j,b_j) \quad \mbox{ with } 0<s_j<r_j<1 \]
such that the rectangles $R_j=f_j(Q)$ are disjoint subsets of $Q.$ The self-affine set generated by $f_1,...,f_m$ is the unique closed nonempty subset $E$ of $Q$ which satisfies
\[  E=f_1(E)\cup ...\cup f_m(E)  ,\]
see \cite{Fal,Bar}. An example is given in Figure 1a. The assumption $s_j<r_j$ will be essential for our results, cf. \cite{BK}, Example 10.

\begin{figure}[ht]
\includegraphics[width=.318\textwidth]{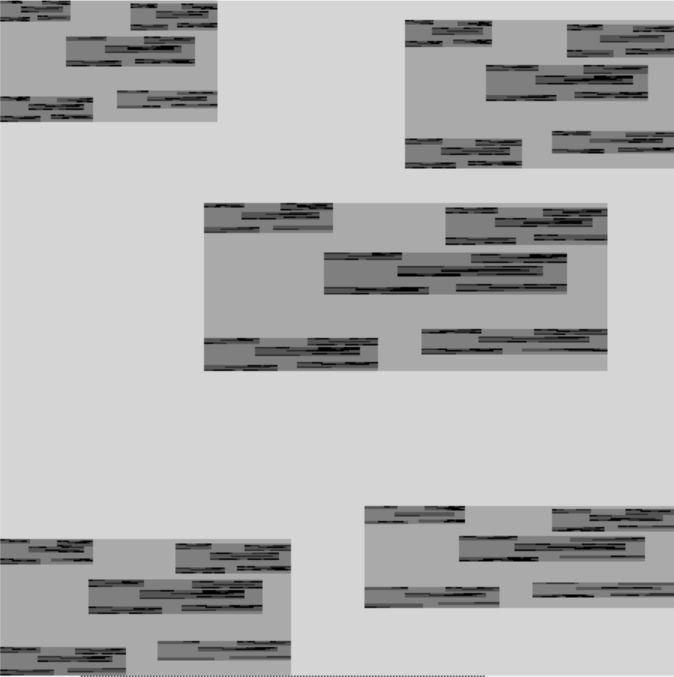}
\includegraphics[width=.32\textwidth]{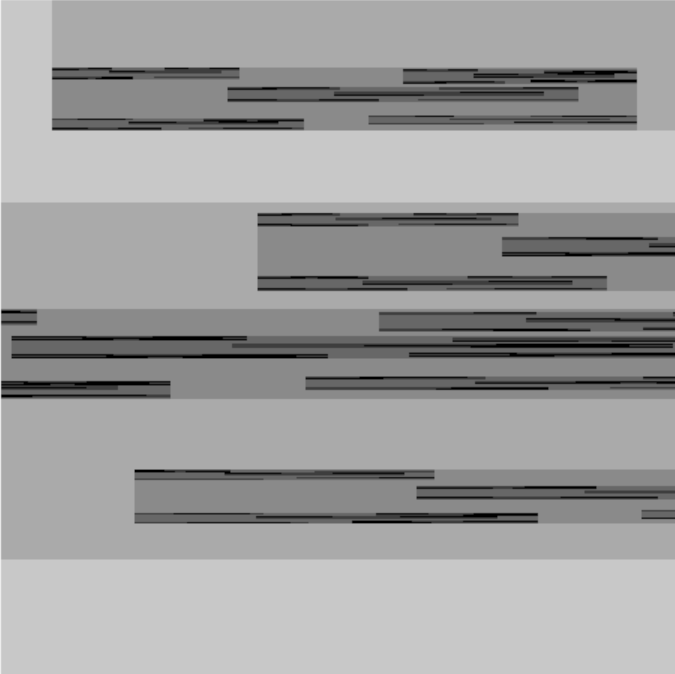}
\includegraphics[width=.32\textwidth]{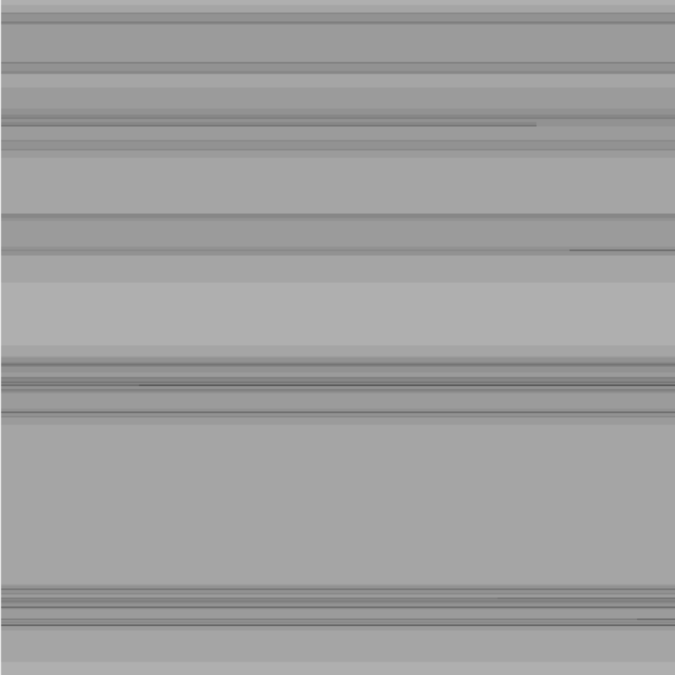}
\vspace{1.5ex}\\
\caption{a) The $R_j=f_j(Q)$ in the unit square $Q.$    b) Magnification, $0.65<x_1<0.77,\, 0.54<x_2<0.66\, .$
c) Magnification of a smaller window $[ 0.73687, 0.736873]\times [0.61145, 0.611453]$ shows the  typical fibre structure which remains preserved under further magnification. }
\end{figure}

For words $u=j_1...j_n\in J^n,$ we consider the mapping $f_u=f_{j_1}\circ ...\circ f_{j_n}$ and the rectangle $R_u=f_u(Q).$ In the figure, the rectangles $R_u$ are shaded in graytones depending only on the length $n$ of the word. 
Thus each greytone colors a set $ E^{(n)}=\bigcup_{u\in J^n} R_u,$ and darker color is used for larger $n.$  The $E^{(n)}$ form a decreasing sequence with $ E=\bigcap_{n=1}^\infty E^{(n)},$ so black color indicates the set $E.$
Since the $R_j$ are disjoint, $E$ is a Cantor set, i.e. totally disconnected without isolated points.  Actually, the address map
$\pi:J^\infty\to E$  defined by  $\pi(j_1,j_2,...)=\bigcap_{n=1}^\infty R_{j_1...j_n}$ 
is a homeomorphism, see \cite{Bar}.

If we magnify $E,$ the local structure appears to be a product set of interval $[0,1]$ and a Cantor set $C\subset \R ,$ as indicated in Figure 1c. Similar computer experiments show that such a local fibre structure exists virtually everywhere.  At first sight, it seems surprising that a Cantor set will develop connected components when it is magnified. Our proof will show that this is due to the different contraction factors in horizontal and vertical direction.  We have chosen simple assumptions to produce a clear argument. We conjecture, however, that the local fibre structure is typical for general self-affine sets where the $f_i$ can involve different rotations and the $R_i$ can heavily overlap.\smallskip

\emph{Tangent sets.} To formulate our result in rigorous terms, we use tangent sets which are defined
like tangential measures \cite[Chapter 14]{Mat}, but in a topological setting. A tangent set is a bit more special than a     micro-set of Furstenberg  \cite{Fu}. The difference is that we fix a basic point $x\in E$ at which magnification will proceed, as in \cite{BFU,MP,Ga,Hoc}. Due to our rectangular setting, it will be convenient to define bounded tangent sets as subsets of the unit square $Q,$ as in \cite{Fu}.  All tangent sets contain the center of $Q$ which corresponds to the basic point $x.$

For a point $x=(x_1,x_2)$ in $E$ and a positive number $t<\frac12$ let $Q_{x,t}$ we consider the square neighborhood 
$Q_{x,t}=[x_1-\frac{t}{2}, x_1+\frac{t}{2}]\times [x_2-\frac{t}{2}, x_2+\frac{t}{2}] $
and the normalizing map \[ h_{x,t}=h: Q_{x,t}\to Q\qquad \mbox{ with }\quad  h(y)=\frac{y-x}{t} +(\textstyle{\frac12,\frac12 })\,  .\]  The set 
\[ N_{x,t}=h(E\cap Q_{x,t})\subset Q\]  is the \emph{normalized view} of $E$ in the neighborhood $Q_{x,t}.$ A \emph{tangent set} of $E$ at $x$ is a limit, with respect to Hausdorff metric (cf. Section 2), of a sequence of normalized views $N_{x,t_n}$ with $t_n\to 0.$  

If $E$ were a differentiable curve, there would be exactly one tangent set at $x,$ namely the line segment through the center $(\frac12,\frac12)$ of $Q$ with the slope of the tangent line of $E$ at $x.$
Due to compactness, tangent sets do exist at all points $x$ of any set $E.$   
In the case of self-similar fractals, there is a large family of tangent sets at $x,$ and these families will coincide for all `typical' points $x.$  All tangent sets of a self-similar fractal Cantor set are isometric to rather large subsets of $E$  (cf. \cite{Ban}).
\smallskip

\emph{Bernoulli measures.} For truly self-affine Cantor sets $E,$ we want to show that tangent sets contain connected fibres: all tangent sets at $x$ have a product structure $[0,1]\times C$ as in Figure 1c.  However, if $x$ is on the left borderline of $Q$ or of some $R_u,$ then all tangent sets at $x$ will contain no point $(y_1,y_2)$ with $y_1<\frac12 .$ Thus we have to confine ourselves to  `typical' points $x.$  

To this end, we consider a Bernoulli measure $\nu_p$ on the address space $J^\infty$ defined by some probability vector $p=(p_1,...,p_m)$ with $p_j>0$ and $\sum p_j=1 .$ That is, $\nu_p\{ (j_1,j_2...)|\  j_k=i_k\mbox{ for } k=1,...,n\} =
p_{i_1}\cdot ...\cdot p_{i_n} .$ Since the address map $\pi :J^\infty\to E$ is one-to-one, we can consider $\nu_p$ as a probability measure for the points $x$ of $E.$  Let $M$ denote the maximal number of rectangles $R_j, j\in J$ which intersect a vertical segment $\{ x_1\}\times [0,1],$ taken over $0\le x_1\le 1.$
Now we can formulate our main result. 

\begin{Thm}
Let $E\subset Q$ be the self-affine Cantor set corresponding to the maps
$ f_j(x_1,x_2)=(r_jx_1, s_jx_2)+(a_j,b_j), j=1,...,m$ where $0<s_j<r_j<1,$ and the rectangles $R_j=f_j(Q)$ are disjoint subsets of $Q.$  We assume that for some ${\tilde{n}}\ge 1,$ each vertical segment $\{ x_1\}\times [0,1]$ intersects at least two rectangles $R_u$ with $u\in J^{\tilde{n}}.$  Moreover, let $\nu_p$ be the Bernoulli measure on $\{ 1,...,m\}^\infty$ associated with an arbitrary probability vector $p=(p_1,...,p_m),$ with $0<p_j<1/M$ for all $j.$

Then for $\nu_p$ almost all points $x$ in $E,$ all tangent sets of $E$ at $x$ have the form $[0,1]\times C$ where $C\subset [0,1]$ is a Cantor set.
\end{Thm}

In Figure 1 we have $M=3,$ and the condition on vertical segments is fulfilled for $\tilde{n}=2.$ This condition requires that two pairs of rectangles touch the right and left border of $Q,$ respectively.  It is used to guarantee that a line segment cannot be a tangent set (cf. Proposition 2). This condition as well as the requirement $p_j<1/M$ can probably be relaxed by using measure-theoretic methods. 
We note that the theorem can be generalized to the $d$-dimensional setting, using maps $ f_j(x_1,x_2,...,x_d)=(r_jx_1, s_{2j}x_2,..., s_{dj}x_d)+(b_{1j},b_{2j},...,b_{dj})$ with $0<s_{kj}<r_j<1.$ The proof of Theorem 1 is given in Section 2.\smallskip
 
The next theorem says that our self-affine sets have ``uniformly scaling scenery'' in the sense of Furstenberg, Gavish, and Hochman.

\begin{Thm}
Under the assumptions of Theorem 1, there is a set $\tilde{E}\subset E$ with $\nu_p( \tilde{E})=1$ for every Bernoulli measure $\nu_p$ with $0<p_j<1/M$  such that\\  for each $x\in \tilde{E},$ each tangent set $T$ at $x$ is a tangent set at each other point $y$ of $\tilde{E}.$\\
Thus the family of all tangent sets does neither depend on the point $x$ nor on the particular measure $\nu_p .$
\end{Thm}

We prove that $E\setminus \tilde{E}$ is a zero set for many Bernoulli measures. One might ask if the Hausdorff dimension of this exceptional set is small. Unfortunately our technique does not provide such
estimates. It is also worth noticing that there exists Bernoulli
measures satisfying the assumptions of Theorem 2 and having arbitrary
small Hausdorff dimension. Although zero sets of such measures can be
large in other senses, we feel that the notion of smallness in Theorem 2
is natural since Bernoulli measures form the geometrically meaningful
family of explicitly known measures for our setting of self-affine sets.

In Section 3 we prove Theorem 2 and discuss the structure of the Cantor sets $C$ in the tangent set $T=[0,1]\times C.$ It will be seen that $C$ is almost a subset of $E.$ 
Under some additional conditions we specify the distribution with respect to $\nu_p$ of the random Cantor set $C.$ 

\section{Structure of tangent sets}
We introduce some notation. Let  \[ \delta =\\ \min\{ |x-y|\, |\, x\in R_i\mbox{ and }y\in R_j\mbox{ for some } i,j\in J \mbox{ with }i\not= j \}\, , \]  
$s^*=\max_{j=1}^m s_j$ and $s_*=\min_{j=1}^m s_j.$ Similarly, $r^*,r_*$ and $p^*,p_*$ denote the maximum and minimum of the $r_j$ and $p_j,$ respectively. We first prove that there is a uniform lower bound for the vertical extension of $E$ in all normalized views of $E.$  

\begin{Prop} Assume the conditions of Theorem 1. 
\begin{enumerate}
\item[(i) ]  For $0\le x_1\le 1$ the set $(\{ x_1\}\times [0,1])\cap E$ has diameter at least $\delta s_*^{\tilde{n}-1} .$
\item[(ii)] For a normalized view $N_{x,t}$ with $t<\delta s_*^{\tilde{n}-1},$ the set $(\{ \frac12\}\times [0,1])\cap N_{x,t}$  has diameter at least $\frac12\delta s_*^{2\tilde{n}-1}.$  
\end{enumerate}
\end{Prop}

\begin{Proof} (i): We use the condition on vertical segments. It implies that the sets $E^{(k\tilde{n})}, k=1,2,...$ intersect each segment $\{ x_1\}\times [0,1]$. So this is true also for $E,$ and by self-affinity for each set $E\cap R_u.$ Thus if the segment intersects both $R_u$ and $R_v,$ with $u,v  \in J^{\tilde{n}},$  it contains a point of $E$ in both of the rectangles. The distance of $R_u$ and $R_v$ on the segment is at least $\delta s_*^n$ where $n$ is the length of the common prefix of $u$ and $v.$  Since $n\le \tilde{n}-1 ,$ the estimate is true.  

(ii): Let $k$ be the largest integer for which $A:=(\{ x_1\}\times [0,1])\cap Q_{x,t}$ intersects only one rectangle $R_u$ with $u\in J^{k\tilde{n}}.$ The assumption $t<\delta s_*^{\tilde{n}-1}$ together with the proof of (i) implies $k\ge 1.$ We apply the proof of (i) again to the two subrectangles of $R_u$ of order $(k+1)\tilde{n}$ which intersect $A.$ We obtain
\[\diam A\cap E \ge s_u\cdot\delta s_*^{\tilde{n}-1}\,  .\]
If (ii) did not hold for $N_{x,t},$ the above diameter would be smaller than $t\cdot \frac12\delta s_*^{2\tilde{n}-1}$ which implies $\frac{t}{2}>\frac{s_u}{s_*^{\tilde{n}}}.$ Now if $u'$ denotes the prefix of $u$ of length $(k-1)\tilde{n},$ we have $\frac{t}{2}>s_{u'}.$ Since $x\in R_{u'},$ this means that $(\{ x_1\}\times [0,1])\cap R_{u'}$ is contained in $Q_{x,t}.$
On the other hand the vertical segment condition says that beside $R_u$ another subrectangle $R_v\subset R_{u'}$ with 
$v\in J^{k\tilde{n}}$ will intersect  $\{ x_1\}\times [0,1].$ This contradicts the choice of $k.$
\end{Proof}

Let us recall the definition of the Hausdorff distance between two sets $A$ and $B$ (cf.~ \cite{Bar}).  We have
$d_H(A,B)\le \varepsilon$  if and only if for every point $x\in A$ there is $y\in B$ with $|x-y|\le \varepsilon,$ and for each  $y\in B$ there is  $x\in A$ with $|x-y|\le \varepsilon .$  \\
For given $\varepsilon>0,$ a subset $S$ of $Q$ is called an \emph{$\varepsilon$-pattern} if $S=[0,1]\times\bigcup_{k=1}^\ell I_k$ where $\ell$ is a positive integer and the $I_k$ are disjoint closed intervals of positive length $\le\varepsilon .$  Proposition 5 below will show that normalized views of $E$ are near to such patterns. The following statement then implies that tangent sets have product structure.

\begin{Prop} Let $\varepsilon_i>0, i=1,2,...$ be a sequence converging to 0, and let $S_i$ be a sequence of $\varepsilon_i$-patterns and $T\subset Q$ the limit set, with $d_H(S_i,T)<\varepsilon_i$ for all $i.$ If the area of $S_i$ converges to zero, $T$ has the form $T=[0,1]\times C$ where $C$ is a nowhere dense subset of $[0,1]$ with zero length. 
\end{Prop}
\begin{Proof} For each $i,$ $T$ is contained in the $\varepsilon_i$-neighborhood $S_i'=[0,1]\times\bigcup_{k=1}^{\ell_i} I_k'$ of $S_i.$  Here $I_k'=[c_k-\varepsilon_i, d_k+\varepsilon_i]$ when $I_k=[c_k,d_k].$ Thus $T\subset\bigcap_{i=1}^\infty S_i'$ which has the form $[0,1]\times C$ with a closed set $C.$ Since the area of $S_i'$ is at most three times the area of $S_i,$ the area of $T$ is zero. So $C$ has length zero and is nowhere dense. 
\end{Proof}

\emph{Definition of approximate normalized views $P_{x,t}^K $.}  Let $N_{x,t}=h(E\cap Q_{x,t})$  with $x\in E$ and $t<\frac12 $ be a normalized view. Let $j_1j_2...$ denote the address of $x,$ and let $n=n(x,t)$ denote the largest integer for which $Q_{x,t}$ intersects only one rectangle $R_u$ with $u=j_1...j_n\in J^n.$ Then $Q_{x,t}$ intersects two different subrectangles $R_{uj_{n+1}}$ and $R_{ui}$  with $i\not= j_{n+1}.$  Since their distance is at least $\delta s_u,$ we have
\begin{equation}\label{diag}  \delta s_*^n\le \delta s_u< t\sqrt{2}\   . \end{equation}
The approximation of $N_{x,t}$ with approximation level $K$ is now defined as
\begin{equation}\label{PK} 
P^K_{x,t}=h(E^{(n+K)}\cap Q_{x,t})=h\left ( Q_{x,t}\cap\bigcup_{w\in J^{n+K} } R_w\, \right ) \,  .
\end{equation}
This is a union of rectangles in $Q,$ as indicated in Figure 1 by one greytone. The following technical statement says that for small $t$ and almost all $x,$  none of the rectangles will end within $Q$ and so  $P_{x,t}^K $ will be a $\varepsilon$-pattern. 
This is the key to the proof of Theorem 1 at the end of this section.

\begin{Prop} Let the conditions of Theorem 1 be fulfilled, and let $\varepsilon_K =(s^*)^K\cdot\frac{\sqrt{2}}{\delta}$ for some positive integer $K.$  The $P_{x,t}^K$  have the following properties. 
\begin{enumerate}
\item[(i) \ ]  Each $P_{x,t}^K$ is a union of closed rectangles of height at most $\varepsilon_K .$
\item[(ii) ]  $P_{x,t}^K$ contains $N_{x,t}$ and $d_H(N_{x,t},P_{x,t}^K)<\varepsilon_K .$
\item[(iii)] If $P_{x,t}^K$ is an $\varepsilon_K$-pattern, it has area at most
$(1-\delta)^{K-\tilde{k}}$ where $\tilde{k}$ denotes the integer part of $\frac{\log\delta /2}{\log s^*}.$
\item[(iv)] The set $B_K=\{ x\in E\, |\,  P_{x,t}^K\  \mbox{ is not an }\varepsilon_K\mbox{-pattern for arbitrary small } t\}$
has measure  $\nu_p(B_K)=0.$
\end{enumerate}
\end{Prop}

\begin{Proof} (i): $P^K_{x,t}$ is a union of closed rectangles in $Q.$  Since each $R_w$ has height $s_w$ and by (\ref{diag})
\[ s_w \le  s_u\cdot (s^*)^K \le t\cdot \frac{\sqrt{2}}{\delta}\cdot (s^*)^K =t\cdot\varepsilon_K\,  ,\] 
we conclude that the rectangles $h(Q_{x,t}\cap R_w)$ have height at most $\varepsilon_K\, .$ 

(ii): By construction, $P^K_{x,t}\supset N_{x,t},$ and for each point $y\in P^K_{x,t}$ there is a point $z\in N_{x,t}$ in the same $h(R_w)$ with $y_1=z_1,$ by the vertical segment condition. So the Hausdorff distance between   $N_{x,t}$ and $P^K_{x,t}$ is bounded by $\varepsilon_K\, .$ 

(iii): We estimate the area of $P_{x,t}^K.$
If this set is an $\varepsilon_K$-pattern, it consists of rectangles which have their endpoints at  $x_1=0$ and $x_1=1.$
The same holds for $P_{x,t}^{K'}=h(E^{(n+K')}\cap Q_{x,t})$ with $0\le K'<K.$  For each rectangle $R_w$ of height $s_w$ with $w\in J^{n+K'},$ there is at least one empty strip of height at least $\delta s_u$ between the subrectangles $R_{wi}, i\in J.$ Thus when we go from  $P_{x,t}^{K'}$ to  $P_{x,t}^{K'+1}$ the area decreases by a factor smaller or equal $1-\delta .$ When we start with area at most $1$ for  $P_{x,t}^0,$ induction gives area at most $(1-\delta)^K$ for  $P_{x,t}^K .$

However, there may be two rectangles at the upper and lower border of $Q_{x,t}$ for which this argument may not work since they are only partially contained in the view. If $s_u\le t,$ we can do the induction with all $R_u$ instead of the part in the square. If $s_u> t,$ we start the induction with $\tilde{k}$ where $2(s^*)^{\tilde{k}}\le \delta ,$ and we estimate the area including the two possible rectangles $R_w, w\in J^{n+\tilde{k}},$ on the boundary, and their subrectangles in the induction which may be outside the square. The definition of $\tilde{k}$ guarantees that those two $R_w$ together have less area than the empty strip of height $\delta s_u> t\delta .$ So we start with area $\le 1$ at $K'=\tilde{k}$ and perform the induction as before.  This completes the proof of (iii).\vspace{1ex}  

(iv): \emph{Borel-Cantelli argument.}  By (\ref{diag}), $n=n(x,t)$ tends to infinity for fixed $x$ and $t\to 0.$ Thus we can write $ P^K_{x,n}$ instead of $ P^K_{x,t},$ and express $B_K$ using $n$ instead of $t.$
\[  B_K=\{ x\in E\, |\,  P_{x,n}^K\  \mbox{ is not an }\varepsilon_K\mbox{-pattern for infinitely many } n\}\, .\]
We shall determine an integer $n_*$ and for all $n\ge n_*$ an exceptional set $A_{n,K}$ such that 
\[  P^K_{x,n}\  \mbox{ is an $\varepsilon_K$-pattern for }\  x\not\in A_{n,K} \] 
and 
\begin{equation}\label{sumank} 
 \sum_{n=n_*}^\infty \nu_p( A_{n,K}) <\infty \ .
\end{equation}
Then
\[  B_K\subset \{ x\in E\, |\, x\in A_{n,K}\  \mbox{ for infinitely many } n\ge n_* \}\, .\]
By the Borel-Cantelli lemma, the right-hand side is a $\nu_p$ null set. Thus $\nu_p(B_K)=0,$ and (iv) is proved.\vspace{0.5ex}

\emph{Construction of $A_{n,K}.$}  
We investigate for which $x$ and $n$ the set $ P^K_{x,n}$ is an $\varepsilon_K$-pattern. Since the height of the rectangles was at most $\varepsilon_K,$ the requirement is that all rectangles reach from $x_1=0$ to $x_1=1.$
This is not always true, see Figure 1 a,b.  The idea is that for larger $n,$ the rectangle $R_u$  becomes longer and narrower, while for fixed $u$ and $K$ only the endpoints $a_w,b_w$ of the  $R_w=[a_w,b_w]\times [c_w,d_w]$ must be avoided by $Q_{x,t},$ where the number $m^K$ of words $w\in J^{n+K}$ with prefix $u$ is constant. Thus for very large $n,$ the set of exceptional points $x$ should be small.

We choose a sufficiently large integer $L$ such that
\begin{equation}\label{defL}   (p^*M)^L\cdot 4m^K < 1  \, .
 \end{equation} 
This $L$ exists since $p^*M<1$ was assumed in Theorem 1. To define $A_{n,K},$ we fix $u\in J^n$ and consider the definition (\ref{PK}) of $P_{x,n}$ as union of the $R_w$ where $w$ is a word of length $n+K$ with prefix $u.$ We say that $v=ui_1...i_L$ is a forbidden word if $R_v$ intersects a vertical segment $\{ c\}\times [0,1]$ where $c$ is one of the $4 m^K$ values $a_w-\frac{t}{2}, a_w+\frac{t}{2}, b_w-\frac{t}{2}, b_w+\frac{t}{2} .$  Each vertical segment intersects at most $M^L$ rectangles $R_v,$ and $\nu_p(R_v)\le\nu_p(R_u)\cdot {p^*}^L$ for each $v.$ (For a Bernoulli measure, $\nu_p(R_v)=\nu_p(R_v\cap E)=\prod_{k=1}^{n+L} p_{v_k}.$) By (\ref{defL}), the union of all $R_v$ over all forbidden words $v$ has $\nu_p$-measure smaller $\nu_p(R_u).$

Now let $A_{n,K}$ be the union of $R_v\cap E$ over all $u\in J^n$ and all forbidden words $v=ui_1...i_L.$ Note that by self-affinity, the set of forbidden words for each $u\in J^n$ is defined by the same suffixes $i_1...i_L.$  Then  $\sum_{u\in J^n} \nu_p(R_u)=1$ and (\ref{defL}) imply $\nu_p(A_{n,K})<1.$  The set $A_{n,K}$ has the following property. When $x\in E\setminus A_{n,K}$ and $t$ is chosen so that $Q_{x,t}$ intersects only one rectangle $R_u$ with $u\in J^n$ but two rectangles of order $n+1$ then $ P^K_{x,t}$ is an $\varepsilon_K$-pattern. \vspace{0.5ex}

\emph{Convergence of the sum.}
To complete the above argument, we show that all rectangles $R_v$ of order $n+L$ have length greater than $t.$ This will work only for sufficiently large $n.$ Let \[ q=\min_{j=1}^m \frac{r_j}{s_j}\,  .\]
By assumption $q>1.$  For the rectangle $R_u$ the length is at least $q^n$  larger than the height. By Proposition 3 (ii),
the height $s_u$ of $R_u$ is at least  $\frac{t}{2} \delta s_*^{2\tilde{n}-1} .$ So the length of $R_u$ is
\[r_u\ge q^{n}s_u \ge q^n \frac{t}{2} \delta s_*^{2\tilde{n}-1}\, . \]
The condition $r_v\ge r_ur_*^L>t$ now is implied by the following restriction on $n$ :
\begin{equation}\label{defn}  q^n\cdot r_*^L> 2\delta^{-1}s_*^{1-2\tilde{n}}  \, .
 \end{equation} 
Let $n_*$ denote the smallest integer $n$ for which this inequality is fulfilled, with $L$ defined by (\ref{defL}). Then the exceptional set $A_{n,K}$ could be properly defined for $n\ge n_*.$ However, to obtain convergence in (\ref{sumank}), we now modify the definition. We define $A_{n,K}$ with $ \nu_p( A_{n,K})< 1$ in the above way only for $n_*\le n< n_*+\ell$ where $\ell$ is the smallest integer with $q^\ell r_*\ge 1.$ For $n_*+\ell \le n< n_*+2\ell$ we define $A_{n,K}$ by forbidden suffixes of length $L+1$ instead of $L.$ The condition (\ref{defn}) will hold and guarantee that the $R_v$ have length larger $t,$ while  (\ref{defL}) modifies to show that  $ \nu_p( A_{n,K})< p^*M .$

In general, we define $A_{n,K}$ for $n_*+k\ell \le n< n_*+(k+1)\ell$  by forbidden suffixes of length $L+k,$ and obtain 
$ \nu_p( A_{n,K})< (p^*M)^k$ from (\ref{defn}). This implies
\[ \sum_{n=n_*}^\infty \nu_p( A_{n,K}) <\frac{\ell}{1-p^*M} <\infty \]
so that the above Borel-Cantelli argument works.
\end {Proof}    

\noindent \emph{Proof of Theorem 1.}  The set $B=\bigcup_{K=1}^\infty B_K$ fulfils $\nu_p(B)=0.$ Take $x$ outside this
set, and take a tangent set $T=\lim_{t_i\to 0} N_{x,t_i}$ of $E$ at $x.$  Since $x$ is not in $B_K,$ we find for every $K$ a number $t_i$ such that $P_{x,t_i}^K$ is an $\varepsilon_K$-pattern and  $d_H(T,N_{x,t_i})\le\varepsilon_K,$ hence $d_H(T,P_{x,t_i}^K)\le 2\varepsilon_K .$ The areas of these patterns converge to zero. Proposition 4, with $2\varepsilon_K$ instead of $\varepsilon_i,$ now shows that $T=[0,1]\times C$ where $C$ is nowhere dense and has length zero.

It remains to show that $C$ cannot contain isolated points and thus is a Cantor set. Let us assume that $C$ contains an isolated point $c,$ with distance $4\alpha>0$ from the rest of $C.$ The following calculation verifies that there are normalized views of $E$ contained in arbitrary thin horizontal strips, which contradicts Proposition 5.

For each given $\beta\in (0,\alpha]$ there is a number $t$ so that the normalized view $N_{x,t}$ satisfies $d_H(T,N_{x,t})<\beta .$ Thus \[ \{ y\in  N_{x,t}\, |\, |y_2-c|<3\alpha\} =  \{ y\in  N_{x,t}\, |\, |y_2-c|<\beta\} \not=\emptyset\ .\]
We take a point $y$ from this set, let $z=h_{x,t}^{-1}(y)$ and $s=t(\alpha +\beta ).$ Then $N_{z,s}\subset [0,1]\times [\frac12-\frac{2\beta}{\alpha +\beta},\frac12+\frac{2\beta}{\alpha +\beta}]\, .$ For small $\beta$ this contradicts Proposition 3 (ii). {\ \hfill $\Box$\medskip}      

\section{Uniformly scaling scenery}
After proving Theorem 1, we can focus on finding the structure of the Cantor set $C\subset [0,1].$ We restrict ourselves to points $x$ outside the exceptional set $B=\bigcup_{K=1}^\infty B_K$ in the proof of Theorem 1.

\begin{Prop}
The set $C$ in each tangent set $T=[0,1]\times C$ at a point $x\in E\setminus B$ is a limit of sets $C_k$ which are normalized intersections of vertical segments with $E.$ 
\end{Prop}

\begin{Proof}  $T$ is a limit of normalized views $N_{x,t_k}.$ We show that $C$ is the  limit of a subsequence of the vertical sections $C_k=N_{x,t_k}\cap (\{\frac12\}\times [0,1]).$ Using Proposition 5 (iv)  we choose the subsequence $(C_K)$ such that $d_H(T, N_{x,t_K})<\varepsilon_K$ and $P_{x,t_K}^K$ is an $\varepsilon_K$-pattern. Then $P_{x,t_K}^K=[0,1]\times C_K.$ Since $d_H(N_{x,t_K},P_{x,t_K}^K)<\varepsilon_K$ by Proposition 5 (ii), the $P_{x,t_K}^K$ converge to $T.$ Hence $C_K$ converges to $C.$ \end{Proof}

Below we shall consider examples where all sets $C$ are themselves normalized vertical sections of $E.$ However, it is not hard to find examples where this is not always the case. The following argument indicates, however, that these can be rather considered as exceptions.  Given $C_K,$ let $V_K$ denote a vertical segment in the unit square such that $C_K$ is a normalized image of $V_K\cap E.$ Then the midpoint of $V_K$ must belong to $E$ since $(\frac12 , \frac12 )$ belongs to each tangent set. If $V_K\subset R_j$ for some rectangle, then we can take the larger homothetic image $f_j^{-1}(V_K)$ instead of $V_K.$ In other words, we can assume that $V_K$ intersects at least two of the rectangles $R_1,...,R_m.$  So $V_K$ and $C_K$ do not differ much in size.

Theorem 1 and the above argument can be roughly summarized as follows: \emph{ essentially all tangent sets are products of the horizontal component $[0,1]$ and a vertical component $C$ which appears in $E$ on a rather large scale.} For the vertical component, the case is almost the same as for self-similar Cantor sets \cite{Ban}.  Now we prove Theorem 2 which says that each of these tangent sets can be found at essentially all points of $E.$ \smallskip 

\noindent \emph{Proof of Theorem 2. Definition of $\tilde{E}.$}  The set $B=\bigcup_{K=1}^\infty B_K$ in the proof of Theorem 1 has measure zero for all Bernoulli measures $\nu_p$ with $0<p_j<1/M .$  We define $\tilde{E}$ as the set of all $x\in E\setminus B$ for which the address contains all words $w$ from the alphabet $J.$ Of course, if all words of arbitrary length are contained in $j_1j_2j_3...$ then each word $w$ is contained infinitely often.
The ergodic theorem for Bernoulli shifts implies that this property holds for $\nu_p$ almost all addresses $j_1j_2j_3...$ Thus $\nu_p(\tilde{E})=1$ for all the measures $\nu_p.$

\emph{Main idea.} Now let $T$ be a tangent set of $E$ at some point $x\in\tilde{E},$ and $y$ another point in $\tilde{E}.$ We want to show that $T$ is a tangent set at $y.$ For an arbitrary $K,$ we choose with Proposition 5 a normalized view $N_{x,t}$ such that $d_H(T,N_{x,t})<\varepsilon_K$ and $P=P_{x,t}^K$ is an $\varepsilon_K$-pattern.  We shall find arbitrary small $s>0$ such that $N_{y,s}$ fulfils $d_H(P,P_{y,s}^K)<\varepsilon_K .$ For small $s$ Proposition 5 implies 
\[ d_H(T, N_{y,s})\le d_H(T,N_{x,t})+d_H(N_{x,t},P)+d_H(P,P_{y,s}^K)+d_H(P_{y,s}^K,N_{y,s})< 4\varepsilon_K .\]
 Performing this construction for every $\varepsilon_K ,$ we see that $T$ is a tangent set at $y$ and Theorem 2 is true.

\emph{Finding  $N_{y,s}.$} We make an assumption on boundary effects which will be verified below. We assume that the rectangles $R_w$ in the definition (\ref{PK}) of the  $\varepsilon_K$-pattern $P=P_{x,t}^K$ do not end exactly at the left, right, upper or lower border of the square $Q_{x,t}.$   When this condition holds, let $\delta<\varepsilon_K$ smaller than the distance of any side of a rectangle $R_w$ to a parallel border of the  square $Q_{x,t}.$ By assumption $\delta>0 ,$ and for $N_{z,t}$ with $|z-x|<\delta$ the approximation  $P_{z,t}^K$ fulfils
$d_H(P_{z,t},P)<\delta .$ Moreover, we find a prefix $w$ of the address of $x$ such that for all $z\in R_w$ we have $|z-x|<\delta .$

Now we use self-affinity: for all $y$ with an address of the form $uwv$ where $u$ is a word and $v$ a sequence and for $s=r_ut$ the approximation $P_{y,s}^K$ must coincide with some $P_{z,t}$ Since $f_u(R_w)=R_{uw},$ the vertical structure of the strips is the same, and the horizontal structure consists of lines anyway.  By the definition of $\tilde{E},$ the address of the given point $y$ contains the word $w$ infinitely many times, so we find  $N_{y,s}$ with arbitrary small $s.$

\emph{Removing boundary effects.} If the assumption on  $N_{x,t}$ is not satisfied, we use the fact that the ``magnification flow  acts on the tangent sets'' \cite{Hoc}. We take a tangent set $T'$ which contains $T$ in the following sense. $T=h(T'\cap [\frac{\alpha}{2}, 1-\frac{\alpha}{2}]^2)$ where $h(z)=\frac{z-m}{1-\alpha}+m$ with $m=(\frac12,\frac12 )$ and some $\alpha>0.$ Since the space of compact subsets of $Q$ with $d_H$ is compact, $T'$ can be taken as the limit
of a subsequence of $N_{x,t_k'}$ where $t_k'=t_k/(1-\alpha)$ and $N_{x,t_k}$ converges to $T.$ When we require that the level $K$ approximation of $N_{x,t_k'}$ is an $\varepsilon_K$-pattern then the rectangles of $P^K$ will extend beyond left and right borders for every $N_{x,\overline{t_k}}$ with $t_k<\overline{t_k}<t_k' .$ Moreover, for each $k$ there is at most a finite number of exceptional values $\overline{t_k}$ between $t_k$ and $t_k'$ for which  the upper or lower border of the square $Q_{x,\overline{t_k}}$ coincides with  the upper or lower border of a rectangle of approximation level $K.$  Thus we find a number $\beta\in (0,\alpha )$ such that $\overline{t_k}=t_k'(1-\beta)$ is not an exception for any $k.$ The above proof now works for $\overline{T}=h(T'\cap [\frac{\beta}{2}, 1-\frac{\beta}{2}]^2)$  with  $h(z)=\frac{z-m}{1-\beta}+m$ instead of $T.$ It shows that $N_{y,\overline{s_k}}$ converges to $\overline{T}$ for some sequence $(\overline{s_k}).$ Then $N_{y,s_k}$ with $s_k=\frac{1-\alpha}{1-\beta}\cdot \overline{s_k}$ converges to $T.$
{\ \hfill $\Box$\medskip}

Attractors of hyperbolic dynamical systems have the local structure $[0,1]\times C,$ where the Cantor set $C$ changes continuously when we run with $x$ along a trajectory. For our self-affine sets, there are no trajectories. When we change $x_1,$ the set $C$ will change its topological structure whenever $x$ passes the boundary of a rectangle. We shall briefly explain why in our case $C$ should rather be interpreted as a random Cantor set. Instead of the family of all tangent sets, we could also have studied the distribution of the sets $T$ or $C$ with respect to $\nu_p,$ cf. \cite{MP,Gr,Ban}. We better restrict ourselves to a simple class of examples. 

\begin{figure}[ht]
\includegraphics[width=.318\textwidth]{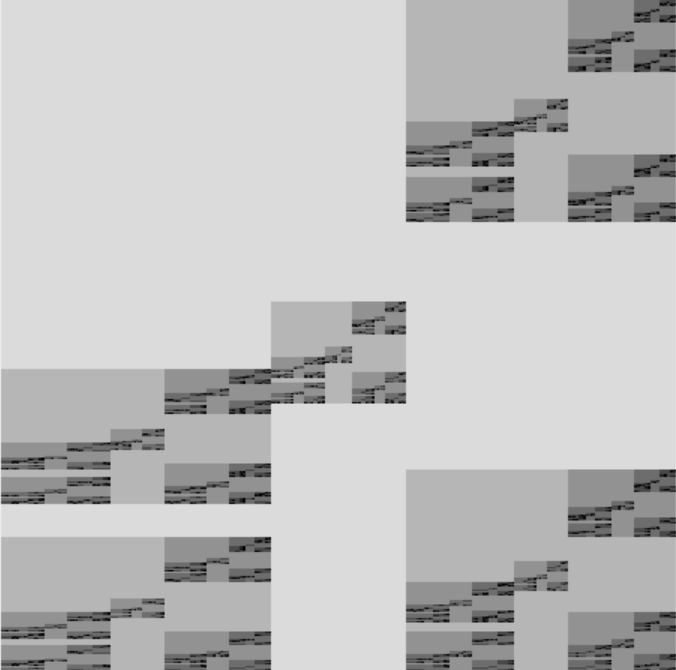}
\includegraphics[width=.32\textwidth]{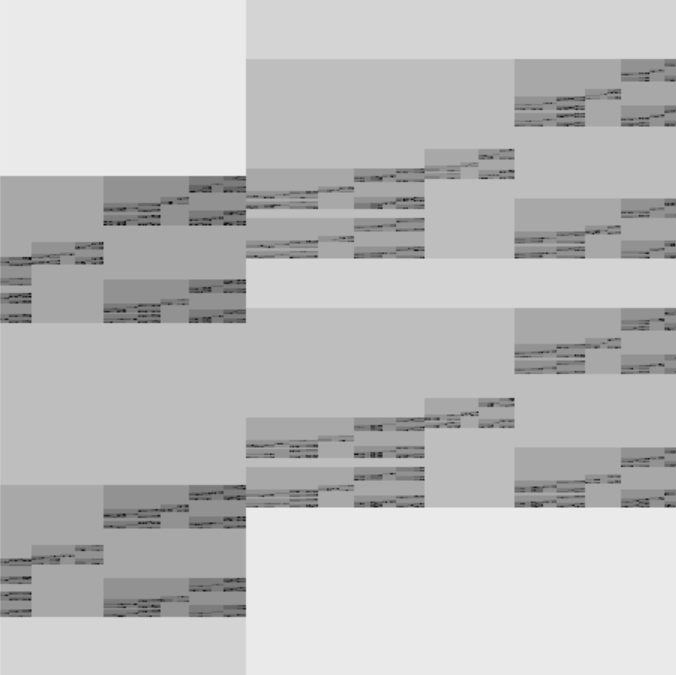}
\includegraphics[width=.32\textwidth]{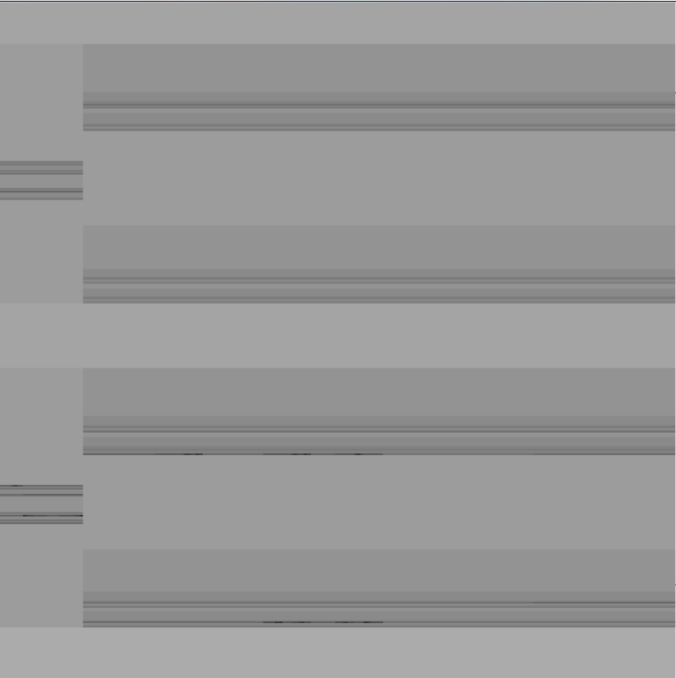}
\vspace{1.5ex}\\
\caption{An example fulfilling the conditions of Proposition 7.  Magnifications of squares with side length $0.12$ and $2\cdot 10^{-8},$ respectively. All vertical sections of $E$ are self-similar random Cantor sets. }
\end{figure}

We assume the following condition due to Gatzouras and Lalley \cite{GL}. Let $\pi(x_1,x_2)=x_1$ denote the projection on the horizontal axis. We assume that the intervals $\pi(R_j), j\in J$ either coincide or have no interior points in common. We also require that the union of the projections is $[0,1],$ which is weaker than the vertical segment condition of Theorem 1.  Figure 2 shows an example.

Thus we have numbers $0=a_1<a_2<...< a_k<a_{k+1}=1$ such that for each $j\in J$ there is an $i_j$ with $\pi (R_j)=[a_{i_j}, a_{i_j+1}].$ Let $J_i=\{ j\in J\, |\, i_j=i\}$  for $i=1,...,k.$ Then each  $f_j$ with $j\in J_i$ has the form
\begin{equation}\label{fj} f_j(x_1,x_2)= (r_ix_1+a_i\, ,\, s_jx_2+b_j)= (h_i(x_1)\, , \, g_j(x_2))\ .
\end{equation}
Here $r_i=a_{i+1}-a_i,$ so the one-dimensional self-similar set generated by $h_1,...,h_k$ is the unit interval. The intervals $h_i([0,1])$ are disjoint, so every point $x_1\in [0,1]$ has an address $i_1i_2... \in \{ 1,...,k\}^\infty$ and only the  endpoints of subintervals have two addresses. Now we can accurately describe the Cantor sets $C.$\vspace{2ex}

\begin{Prop} (cf. \cite{Be})  Let the self-affine set $E$ be generated by the mappings $f_j$ defined in (\ref{fj}). Then for each $x_1\in [0,1]$ with address $i_1i_2... \in \{ 1,...,k\}^\infty ,$ the set \\ $C=E\cap (\{x_1\}\times [0,1])$
has the representation
\[ C=\{ x_2 \,|\, \mbox{ there are } j_k\in J_{i_k},\, k=1,2,...\mbox{ with }  \lim_{n\to\infty} g_{j_1}\circ ...\circ g_{j_n}(0)=x_2 \} \, .\] \end{Prop}

In other words, the sets $J_i, i=1,...,k$ define $k$ families of one-dimensional functions $g_j(t)=s_jt+b_j, \, j\in J_i.$ The set $C$ is the random self-similar subset of $[0,1]$ obtained by applying a function in $J_{i_1}$ to a function in $J_{i_2}$ to a function in $J_{i_3}$ etc. When $x_1$ is one of the countably many points with two addresses, we obtain the union of the two corresponding Cantor sets.  The proof is simple calculation.

To speak of a random construction, we still need a probability measure. For a Bernoulli measure $\nu_p$ on $J^\infty$ with probability vector $(p_1,...,p_m)$ we define $P_i=\sum_{j\in J_i} p_j$  for $i=1,...,k.$ Then at each level, we choose an $i\in\{ 1,...,k\}$ with probability $P_i,$ and apply all the mappings $f_j$ with $j\in J_i.$ This is equivalent to saying that the probability space is $\Omega=[0,1],$ and the probability measure is the self-similar measure generated on $[0,1]$ by the similarity maps $h_i$ and weights $P_i$ with $i=1,...,k.$  A good choice of the $p_j$ depends on the size of the $R_j$ but we do not go into details here.

Fine dimension estimates for random fractal constructions are known to involve logarithmic corrections \cite{GMW}. In fact, the detailed Hausdorff dimension calculations of self-affine Gatzouras-Lalley carpets by Peres \cite{P} involve logarithmic corrections, different from those in  \cite{GMW} due to the different probabilistic construction.

\end{document}